[1994/12/01]
\documentclass[draft]{article}
\title{AF-equivalence relations and their cocycles}
\author{Jean Renault}
\date{8 November 2001}

\usepackage{amsmath,amsthm}

\chardef\bslash=`\\ 

\newtheorem{thm}{Theorem}[section]
\newtheorem{cor}[thm]{Corollary}
\newtheorem{lem}[thm]{Lemma}
\newtheorem{prop}[thm]{Proposition}

\theoremstyle{definition}
\newtheorem{defn}{Definition}[section]

\theoremstyle{remark}
\newtheorem{rem}{Remark}[section]

\newcommand{\eval}[2][\right]{\relax
  \ifx#1\right\relax \left.\fi#2#1\rvert}

\begin{document}
\maketitle
\markboth{AF-equivalence relations and their cocycles}
{AF-equivalence relations and their cocycles}
\begin{abstract}
After a review of some the main results
about hyperfinite equivalence relations
and their cocycles in the measured setting,
we give a definition of a topological
AF equivalence relation. We show that every
cocycle is cohomologous to a quasi-product cocycle.
We then study the
problem of determining the quasi-invariant
probability measures admitting a given
cocycle as their Radon-Nikodym
derivative.\footnote{

Hyperfinite equivalence relations. AF equivalence relations. Cocycles.
Markov measures. Radon-Nikodym derivative. Perron-Frobenius-Ruelle theorem.\ }

\end{abstract} 
\renewcommand{\sectionmark}[1]{}

\section{Introduction.}

Our motivation is a well-known result about some KMS states of the Cuntz algebras.
Recall that the Cuntz algebra ${\cal O}_d$ (here
$d$ is a finite integer at least equal to 2) is the C$^*$-algebra generated by
$d$ elements
$S_1,\ldots,S_d$ satisfying:
\[\begin{array}{cc}
S_j^*S_i&=\delta_{ij}1\\
\sum_{i=1}^d S_iS_i^*&=1.
\end{array}\]
The gauge automorphism group is the one-parameter automorphism group $\alpha_t$ of
${\cal O}_d$ defined by
$$\alpha_t(S_i)=e^{it}S_i.$$
Then $\alpha$ has a unique KMS state and
it occurs at inverse temperature $\beta=\log d$ (see \cite {op:KMS}).

This result can be expressed in the framework of topological dynamics. The Cuntz
algebra ${\cal O}_d$ is the C$^*$-algebra of the one-sided full shift $T_d$ on the
product space $X_d=\prod_1^\infty\{1,\ldots,d\}$. For example, one can use the
groupoid construction for
\,${\cal O}_d$\, given in the Section~III.2 of
\cite{ren:approach} (see also\cite{ren:cuntzlike}). With a view on more general
graph algebras, we realize $X=X_d$ as the infinite path space of the oriented graph
consisting of a single vertex and
$d$ edges, labelled by
$1,\ldots,d$. Thus an element of $X$ is an infinite path $x=x_1x_2\ldots\quad $.
The one-sided shift $T=T_d$ acts on $X$ according to:
$$Tx_1x_2\ldots=x_2x_3\ldots $$
Then
$$G(X,T)=\{(x,m-n,y)\in X\times{\bf Z}\times X:m,n\in{\bf N}, T^mx=T^ny\}$$
has a natural structure of \'etale groupoid  and 
${\cal O}_d=C^*(G(X,T))=C^*(X,T)$. Note that the equivalence relation defined by
$G(X,T)$ on $X$, namely $x$ and $y$ are equivalent if there exist $m,n$ such that
$T^mx=T^ny$ is the ``tail equivalence with lag''.
\smallskip 
Moreover, the gauge automorphism group is the
dual action (see Section~II.5 of \cite{ren:approach})  of ${\bf
T}=\hat{\bf Z}$ on
$C^*(X,T)$ with respect to 
the fundamental cocycle:
\[\begin{array}{cc}
c:G(X,T)&\rightarrow{\bf Z}\\
(x,k,y)&\mapsto k.
\end{array}\] 
Indeed, we have:
$$\alpha_t(f)(\gamma)=e^{itc(\gamma)}f(\gamma)$$
for $f\in C_c(G(X,T)),\ \gamma\in G(X,T),\ t\in{\bf R}$.
\vskip .5cm
A probability measure $\mu$ on $X=G^{(0)}$ defines a state $\phi_\mu$ on
$C^*(G(X,T))$ such that, for $f\in C_c(G(X,T))$,
$$\phi_\mu(f)=\int f_{|G^{(0)}} d\mu.$$
It is known (\cite{ren:approach}, Proposition 5.4 in Section II, 
see also \cite{rue:hyperbolic}, Theorem 2.1), that $\phi_\mu$ is a KMS state for
$\alpha$ at inverse temperature $\beta$ if and only if $\mu$ is quasi-invariant
for $G(X,T)$ with Radon-Nikodym derivative:
$$D_\mu=e^{-\beta c}.\eqno(RN)$$
Moreover, under an assumption of essential freeness which is satisfied here, all
KMS states are of the above form. In this context, the counterpart of the KMS
problem is:
\vskip .3cm
{\bf The  Radon-Nikodym problem .} {\it Given an \'etale locally compact
groupoid $G$ and a cocycle $D\in Z^1(G,{\bf R}_+^*)$, find all quasi-invariant
probability measures $\mu$ which admit $D$ as their Radon-Nikodym derivative.}
\vskip .5cm
When $G=G(X,T)$ is the groupoid of the one-sided full shift $(X_d,T_d)$ and $c$ is
the fundamental cocycle, the equation $(RN)$ has only a solution for $\beta=\log d$
and then $\mu$ is unique. Note that $\log d$ is the topological entropy of the
shift
$T_d$. One can deduce from the Perron-Frobenius-Ruelle theorem a more general
result: under suitable assumptions on the dynamics $(X,T)$ and on
$\varphi\in C(X,{\bf R})$, the cocycle $e^{c_\varphi}$ defined by $\varphi$ is
a Radon-Nikodym derivative if and only if the topological pressure $p(\varphi)$ of
$\varphi$ is zero. The above case corresponds to a constant function
$\varphi$.
\vskip .5cm
The Cuntz algebra ${\cal O}_d$ has a privileged AF subalgebra, namely the
fixed-point algebra ${\cal O}^\alpha_d$ of the gauge automorphism group, which is
uniformly hyperfinite of type $d^\infty$. ${\cal O}_d$ is the crossed product of
${\cal O}^\alpha_d$ by the shift endomorphism. Correspondingly, the groupoid
$G(X,T)$ has a privileged AF subgroupoid, the kernel
$R(X,T)=c^{-1}(0)$ of the fundamental cocycle. It is the (graph of the) ``tail 
equivalence (without lag)'' relation on the infinite path space $X_d$. By
definition, inductive limit techniques are available for AF equivalence relations.
A possible strategy to study the Radon-Nikodym problem on $G(X,T)$ consists in two steps: the
Radon-Nikodym problem for $R(X,T)$ and the extension to $G(X,T)$ of cocycles on $R(X,T)$.

Before presenting these two steps, I shall review some basic facts on hyperfinite
measured equivalence relations which can serve as guidelines to the topological
case.

The reader should be aware that the sole originality of this talk may be the
emphasize on the topological groupoids $G(X,T)$ and $R(X,T)$. These groupoids
have been introduced and used in the measure=theoretical
setting (see e.g. \cite{av:endo} or \cite{kv:walk}).  The KMS problem for Cuntz
and Cuntz-Krieger algebras and some gauge automorphism groups was studied by
D.Evans in
\cite{eva:on},
\cite{eva:on+},
and in particular in
\cite{eva:markov}, from a rather similar point of view. See also \cite{efw:kms}
and \cite{el:KMS}. One should also mention that the Radon-Nikodym problem for
cocycles on transformation groups was studied by K. Schmidt in
\cite{sch:cocycles} but we do not compare here his results with ours.

\section{Hyperfinite equivalence relations.}

We recall here three fundamental results on hyperfinite equivalence relations.
Besides their relevance in the study of hyperfinite factors, they stimulate
further work on Borel and topological equivalence relations.

\subsection{Amenable and hyperfinite equivalence relations.}

The basic reference on measured equivalence relations is \cite{fm:relations}. I
recall briefly the setting. We assume that $(X,{\cal B},\mu)$ is a
standard measured space and that
$R$ is an equivalence relation on $X$ such that:

$(a)$ its classes are countable,

$(b)$ its graph $R$ is a Borel subset of $X\times X$,

$(c)$ the measure $\mu$ is quasi-invariant under $R$.

The last condition means that the measures $r^*\mu$ and $s^*\mu$ on $R$
are equivalent (where
$r,s$ denote respectively the first and the second projections onto $X$ and
$r^*\mu(f)=\int\sum_{y=x}f(x,y)d\mu(x)$). Then, one defines the Radon-Nikodym
derivative of $\mu$ as $D_\mu=r^*\mu/s^*\mu$. The measure $\mu$ is called
invariant when $D_\mu\equiv 1$. The Radon-Nikodym derivative $D_\mu:R\rightarrow
{\bf R}_+^*$ satisfies the cocycle identity:
$$D_\mu(x,y)D_\mu(y,z)=D_\mu(x,z)\quad\hbox{for a.e.} (x,y,z)\in R^{(2)}.$$

\begin{defn}
An equivalence relation $R$ as above is called {\it finite} if its classes are
finite. It is called {\it hyperfinite} if there exists an increasing sequence of
finite equivalence relations $(R_n)$ such that $R=\cup R_n$ up to a nullset.
\end{defn}

Before stating the first theorem, we recall one of the equivalent definitions of
the amenability of a measured equivalence relation.

\begin{defn}
An equivalence relation $R$ as above is called {\it amenable} if there exists a
family
$\{m^x\}_{x\in X}$, where $m^x$ is a state on $l^\infty[x]$ ($[x]$ denotes the
class of $x$), which is
\begin{enumerate}
\item invariant: $(x,y)\in R\Rightarrow m^x=m^y$,
\item measurable: $f\in L^\infty(R)\Rightarrow x\mapsto m^x(f)$ is measurable.
\end{enumerate}
\end{defn}

\begin{thm}(Connes-Feldman-Weiss, \cite {cfw:amenable}) Let $R$ be an equivalence
relation as above. Then it is amenable if and only if it is hyperfinite.
\end{thm}

\subsection{ Classification of hyperfinite equivalence relations.}

Let us review briefly the classification of hyperfinite equivalence
relations. It is a classification up to isomorphism, where an {\it isomorphism}
from
$R_1$ on
$(X_1,{\cal B}_1,\mu_1)$ onto $R_2$ on $(X_2,{\cal
B}_2,\mu_2)$ is a Borel isomorphism $\phi$ from a conull subset of $X_1$ onto a
conull subset of $X_2$ such that $\phi_*\mu_1\sim\mu_2$ and
$\phi\times\phi(R_1)=R_2$ up to nullsets. A measured equivalence relation $R$ over
$(X,\mu)$ is called {\it ergodic} if the real valued
Borel functions on
$X$ such that
$f\circ r=f\circ s\  a.e$ are constant
$a.e.$   Since an arbitrary measured equivalence relation
can be decomposed into ergodic components over its standard quotient, it suffices
to consider ergodic measured equivalence relations. Ergodic measured equivalence
relations are classified first into
\begin{enumerate}
\item type $I_n$ $(n=1,2,\ldots,\infty)$: $X$ has cardinality $n$ and $R$ is
transitive,
\item type $II_n$ $(n=1\  \hbox{or}\  \infty)$: $R$ is not transitive and there
exists an invariant measure (finite or infinite) equivalent to $\mu$.
\item type $III$: there is no
invariant measure  equivalent to $\mu$.
\end{enumerate}
The type I relations are clearly hyperfinite. Concerning the other types, we have:

\begin{thm} Classification of hyperfinite ergodic equivalence relations.
\begin{enumerate}
\item (Dye \cite{dye:I,dye:II}) Up to isomorphism, there is one and only
one  hyperfinite ergodic equivalence relation of type $II_1$;
\item (Krieger \cite {kri:flows}) Type $II_\infty$ and type $III$
hyperfinite ergodic equivalence relations are completely classified by an
associated flow, namely the Mackey range of the Radon-Nikodym derivative $D_\mu$,
and any ergodic flow can occur.
\end{enumerate}
\end{thm}

Let us recall the construction of the Mackey range of a cocycle. We are given
a measured equivalence relation $(X,{\cal B}, R,\mu)$, a locally compact
group  $A$ with Haar measure $\lambda$ and a Borel
cocycle $c:R\rightarrow A$. One defines
the skew product equivalence relation $R(c)$ on $(A\times X,\lambda\times\mu)$ by
$$(ac(x,y),y)\sim (a,x)\quad\hbox{for}\ (x,y)\in R, a\in A.$$
The measured equivalence relation $R(c)$ is not necessarily ergodic, even if
$R$ is so. The ergodic decomposition of an arbitrary measured equivalence relation
$(X,R,\mu)$ provides the {\it standard quotient} $X//R$. It is a measured space
such that:
$$L^\infty(X//R)=\{f\in L^\infty(X): f(x)=f(y)\quad\hbox{for a.e.} (x,y)\in
R\}.$$ In our case, the standard quotient
$P(c)=A\times X//R(c)$ is a measured space equipped with a left action of $A$.
\begin{defn}
The measured $A$-space $P(c)$ is called the {\it Mackey range} of the cocycle $c$.
\end{defn}
In the case of the Radon-Nikodym cocycle $D_\mu$, it is customary to consider the
cocycle
$c=\log D_\mu:R\rightarrow {\bf R}$; its Mackey range is also called its
{\it Poincar\'e flow}.

An important example of a Mackey range is the Poisson boundary of a random walk on
a group (\cite {zim:poisson}). Recall that a random walk is determined by a
probability measure $m$ on a group $A$. We assume here that the group $A$ is
discrete countable and that the measure $m$ has a finite symmetric support $S$. We
consider the tail equivalence relation $R$ on the measure space $(X=\prod_1^\infty
S,\mu=\prod_1^\infty m)$. It is a hyperfinite ergodic equivalence relation. The
random walk cocycle
$c:R\rightarrow A$ is defined by
$$c(x,y)=\lim_{n\to\infty} x_1x_2\ldots x_ny_n^{-1}\ldots y_2^{-1}y_1^{-1}$$
(the limit exists because the sequence is stationary). The measured $A$-space
$P(c)$ is the Poisson boundary of the random walk $(A,m)$. It describes the bounded
$m$-harmonic functions on $A$. It is worth noting that the tail equivalence
relation $R$ is the tail equivalence relation of the graph having a single vertex
and
$S$ as its set of edges while the skew product while the skew product equivalence
relation $R(c)$ is the
tail equivalence relation of the Cayley graph of $(A,S)$.

\subsection {Quasi-product cocycles.}

The random walk cocycle defined above is an example of a quasi-product cocycle, which we
define now. The definitions given below make sense for an arbitrary locally finite
oriented graph; we restrict ourselves to the case of a Bratteli diagram.
Recall that a Bratteli diagram is an oriented graph $E\rightarrow V$ where the
vertices are stacked on levels $n=0,1,2,\ldots$ and the edges run from a vertex of
level $n-1$ to a vertex of level $n$. We denote by $V(n)$ the set of vertices of
the level $n$ and by $E(n)$ the set of edges from level $n-1$ to level $n$. We
assume that for every vertex
$v$, there are finitely many, but at least one, edges starting from $v$. An
infinite path is a sequence of connected edges
$x=x_1x_2\ldots$, where $x_1$ starts at level $0$. The space $X$ of infinite paths
has a natural totally disconnected topology, with the cylinder sets
$Z(x_1x_2\ldots x_n)$ as a basis; we equip $X$ with the underlying Borel
structure
$\cal B$. The tail equivalence relation $R$ has a Borel graph in
$X\times X$ (this will be made clear in the next section).

\begin{defn}
A cocycle $c:R\rightarrow A$, where $A$ is a group, is called a
{\it quasi-product cocycle} (relative to the Bratteli diagram $(V,E)$) if there exists a
map
$f:E\rightarrow A$, often called a labelling, such that
$$c(x,y)=\lim_{n\to\infty} f_1(x_1)f_2(x_2)\ldots f_n(x_n)[f_n(y_n)]^{-1}\ldots
[f_2(y_2)]^{-1}[f_1(y_1)]^{-1}$$
where $f_n$ is the restriction of $f$ to $E(n)$.
\end{defn}

\begin{defn}
A probability measure $\mu$ on $X$ is
called a {\it Markov measure} (relative to the Bratteli diagram $(V,E)$) if there
exists a map
$p:E\rightarrow {\bf R}_+^*$, called a transition probability, satisfying
$\sum_{s(e)=v}p(e)=1$ for every $v\in V$ and a map $\mu_0:V(0)\rightarrow {\bf
R}_+^*$, called an initial probability, satisfying $\sum_{V(0)}\mu_0(v)=1$ such
that
$$\mu(Z(x_1x_2\ldots x_n))=\mu_0\circ s(x_1)p_1(x_1)p_2(x_2)\ldots p_n(x_n).$$
\end{defn}

Note that a Markov measure $\mu$ is quasi-invariant with respect to the tail
equivalence relation $R$ (its Radon-Nikodym derivative will be given in the
next section). Thus, it turns $R$ into a measured equivalence relation. Together
with a quasi-product cocycle, this provides the most general model for a Borel cocycle on
a hyperfinite equivalence relation:

\begin{thm}(\cite{eg:amenable,aeg:amenable}) Given a hyperfinite equivalence
relation
$(X,R,\mu)$, a locally compact group $A$ and a Borel cocycle $c:R\rightarrow A$,
there exists a Bratteli diagram $(V,E)$ with infinite path space $\underline X$
and tail equivalence relation $\underline R$, a Markov measure $\underline\mu$ on
$\underline X$ and an isomorphism $\phi$ from $R$ onto $\underline R$ carrying $c$
into a cocycle cohomologous to a
quasi-product cocycle $\underline c:\underline R\rightarrow A$ .
\end{thm}

\section{AF equivalence relations.}

As it is well-known, the theories of hyperfinite factors and of hyperfinite
equivalence ergodic equivalence relations are intimately related. The
C$^*$-algebraic analog of a hyperfinite von Neumann algebra is an
approximately finite dimensional (or AF) C$^*$-algebra. There is a notion of
topological AF equivalence relation which is the analog of a hyperfinite measured
equivalence relation. Again, the theories of AF C$^*$-algebras and of AF
equivalence relations are intimately related. In particular, they share the same
complete invariant, namely the dimension range.

\subsection{AP and AF equivalence relations.}

In this section, $X$ is a locally compact second countable Hausdorff space and $R$
is an equivalence relation on $X$. For the sake of simplicity, we assume that the
equivalence classes are countable. As before, we also denote by $R\subset X\times
X$ its graph.

\begin{defn} The equivalence relation $R$ on $X$ will be called:
\begin{enumerate}
\item {\it proper} if its graph $R$ is closed in $X\times X$ and the quotient map
$\pi:X\rightarrow X/R$  is a local homeomorphism;
\item {\it approximately proper} (or AP) if $R=\cup R_n$, where $(R_n)$ is an
increasing sequence of proper equivalence relations;
\item {\it approximately finite} (or AF) if it is AP and $X$ is totally
disconnected.
\end{enumerate}
\end{defn}

\begin{rem}
\begin{enumerate}
\item The inductive limit topology turns an AP equivalence relation $R$ into an
\'etale locally compact Hausdorff groupoid, which is topologically amenable.
Unless we specify otherwise, we equip $R$ with this topology.
\item There may exist distinct topologies on an equivalence relation $R$ making it
into an \'etale locally compact groupoid. Examples can be deduced from
\cite{gps:orbitequivalence}, Theorem 2.3. There
$X$ is the Cantor space and $R$ is the orbit equivalence relation of a minimal
homeomorphism $T$ of $X$ into itself. The authors show that $R$ is AF but its AF
topology is distinct from the topology of $X\times {\bf Z}$.
\end{enumerate}
\end{rem}

The tail equivalence relation on the infinite path space of a Bratteli diagram is
an example of an AF equivalence relation and a quasi-product cocycle an example of
a continuous cocycle. In fact, just as in 2.3 above, these examples exhaust all AF
equivalence relations and their continuous cocycles (more precisely, their
cohomology class). The proof is similar to that of \cite{eg:amenable,aeg:amenable}.

\begin{thm} Let $(X,R)$ be an AF equivalence relation  with $X$ compact and
let $A$ be a locally compact group. 
\begin{enumerate}
\item There exists a Bratteli diagram with tail equivalence relation $(\underline
X,\underline R)$ and a homeomorphism $\varphi:
X\rightarrow\underline X$ carrying $R$ into $\underline R$.
\item If moreover $c\in Z^1(R,A)$ is
a continuous cocycle taking finitely many values on each
$R_n$ of an increasing sequence of proper relations such that $R=\cup R_n$, the
Bratteli diagram can be chosen such that $\varphi$ carries
$c$ into a quasi-product cocycle.
\item Every continuous cocycle $c\in Z^1(R,A)$ is continuously cohomologous to
a continuous cocycle taking finitely many values on each
$R_n$ of a given increasing sequence of proper relations such that $R=\cup R_n$
\end{enumerate}
\end{thm}
\begin{proof}
Given a surjective map $\pi:X\rightarrow Y$,
we say that a subset $\omega subset X$ is a $\pi$-section if the restriction of
$\pi$ to $\omega$ is injective. We say that a partition $\cal P$ of $X$ by
$\pi$-sections is a
$\pi$-partition if for all $\omega,\omega'\in {\cal P}$, either
$\pi(\omega)\cap \pi(\omega')=\emptyset$ or $\pi(\omega)=\pi(\omega')$. We denote
by
$\pi_*{\cal P}$ the partition $\{\pi(\omega), \omega\in{\cal P}\}$ of $Y$. If $X$
is compact and totally disconnected and $\pi$ is a local homeomorphism, for every
partition
$\cal Q$ of
$X$ by clopen sets, there exists a $\pi$-partition of $X$ which is finer than
$\cal Q$. 

Let $(X,R)$ be an AF equivalence relation. Thus $R=\cup R_n$ where
$(R_n)$ is an increasing sequence of proper equivalence relations on $X$ (we
assume that $R_0$ is the diagonal) . We write
$X_n=X/R_n$  and
we denote by $\pi_n:X\rightarrow X_n$ the quotient maps. We also have maps
$\pi_{n,n-1}:X_{n-1}\rightarrow X_n$ such that
$\pi_n=\pi_{n,n-1}\circ\pi_{n-1}$. Then, there exists for all $n\ge 1$ a
$\pi_{n,n-1}$-partition ${\cal P}_n$ of $X_{n-1}$ by clopen sections such
that ${\cal P}_n$ refines $(\pi_{n-1,n-2})_*{\cal P}_{n-1}$. Moreover, if
partitions  ${\cal Q}_n$ of $X_{n-1}$ by clopen subsets have been specified,
one can choose $({\cal P}_n)$ such that ${\cal P}_n$ refines ${\cal Q}_n$. Such a
sequence of partitions $({\cal P}_n)$ wil be called a tower. A tower defines a
Bratteli diagram. Namely, we define
$V(0)$ as single vertex and for $n\ge 1$,
$V(n)=(\pi_{n-1,n-2})_*{\cal P}_{n-1}$, $E(n)={\cal P}_n$ and for
$\omega\in{\cal P}_n$, $s(\omega)$ is the element of $(\pi_{n-1,n-2})_*{\cal
P}_{n-1}$ containing $\omega$ while $r(\omega)=\pi_{n,n-1}(\omega)$. Let
$(\underline X,\underline R)$ be the AF equivalence relation associated to this
Bratteli diagram. An element of
$\underline X$ is a sequence
$(\omega_n)$, where $\omega_n$ is a clopen subset of
$X_{n-1}$, such that $\omega_{n+1}$ is contained in $\pi_{n,n-1}(\omega_n)$. The
map
$\varphi:X\rightarrow \underline X$ associate to $x\in X$ the sequence
$(\omega_n)$, where $\omega_n$ is the element of the partition ${\cal P}_{n-1}$
which contains $\pi_{n-1}(x)$. If $\{\pi_n^{-1}({\cal
P}_n)\}, n\in{\bf N}$ is a basis for the topology of $X$, then $\varphi$ induces
an isomorphism of the AF equivalence relations $(X,R)$ and $(\underline
X,\underline R)$.

Suppose now that
$c:R\rightarrow A$ is a continuous cocycle taking finitely many values on $R_n$
for all $n$ (note that this necessarily the case when $G$ is discrete since $R_n$
is compact). We choose continuous sections
$\sigma_{n,n+1}$ for the maps
$\pi_{n,n+1}$ (by specifying sections for $r:E(n)\rightarrow V(n)$). By
composition, we get continuous sections $\sigma_n$ for $\pi_n:X\rightarrow X_n$.
We define $b_n:X\rightarrow A$ by $b_n(x)=c(x,\sigma_n\circ\pi_n(x))$ and 
$b'_n:X_n\rightarrow A$ by
$b'_n(x_n)=c_n(x_n,\sigma_{n,n+1}\circ\pi_{n,n+1}(x_n))$, $x_n\in X_n$ and
$c_n(x_n,y_n)=c(\sigma_n(x_n),\sigma_n(y_n))$. Then
$c(x,y)=b_n(x)b_n(y)^{-1}$ for all $(x,y)\in R_n$ and
$b_{n+1}(x)=b_n(x)b'_n\circ\pi_n(x)$. Since $b'_n$ takes only
finitely many values, we can construct by induction a tower ${\cal Q}_n$ such that
$b'_n$ is constant on the elements of the partition ${\cal Q}_n$. Then, $c$ is a
quasi-product cocycle relative to the Bratteli diagram associated with this tower.

The last part of the theorem is obtained by successive applications of the lemma
below. We can assume that $A$ is metrizable and endowed with a complete right
invariant metric
$d$. We specify a sequence of positive numbers
$(\epsilon_n)$ such that
$\sum\epsilon_n<\infty$ and construct by induction continuous functions
$\varphi_n: X\rightarrow A$ such that $d(\varphi_n(x),1)\le\epsilon_n$ for all
$x\in X$ and 
$$c_n=(\varphi_n\circ r)\ldots(\varphi_1\circ r)c(\varphi_1\circ s)^{-1}\ldots
(\varphi_n\circ s)^{-1}$$
takes finitely many values on $R_n$ and agrees with $c_{n-1}$ on $R_{n-1}$.
Then\hfill\break
$\varphi_n(x)\ldots \varphi_1(x)$ converges uniformly to a limit $\psi(x)$ and
$(\psi\circ r)c(\psi\circ s)^{-1}$ takes
finitely many values on each $R_n$.

\end{proof}

\begin{lem} Let $X,Y,Z$ be totally disconnected compact spaces,
$\pi:X\rightarrow Y$ and $\rho:Y\rightarrow Z$ surjective local homeomorphisms,
$S=X\times_Y X$ and $T=X\times_Z X$ and let $G$ be a group endowed with a right
invariant metric $d$. Given a continuous cocycle $c\in Z^1(T,A)$
taking finitely many values on $S$ and $\epsilon>0$, there exists $\varphi:
X\rightarrow A$ continuous such that $d(\varphi(x),1)\le\epsilon$ for all $x\in X$
and
$(\varphi\circ r)c(\varphi\circ s)^{-1}$ agrees with $c$ on $S$ and takes
finitely many values on
$T$.

\end{lem}
\begin{proof} We fix continuous sections $\sigma$ for $\pi$ and $\tau$ for $\rho$.
We define $b:X\rightarrow A$ by $b(x)=c(x,\sigma\circ\pi(x))$ and $a:Y\rightarrow
A$ by $a(y)=c(\sigma(y),\sigma\circ\tau\circ\rho(y))$ .
There exists $a':Y\rightarrow A$ continuous and finitely valued such that
$d(ta'(y)(a(y))^{-1}t^{-1},1)<\epsilon$ for all $y\in Y$ and all $t$ in the range
of
$b$. Then the function $\varphi:X\rightarrow A$ defined by 
$$\varphi(x)=b(x)[a'\circ\pi(x)(a\circ\pi(x))^{-1}](b(x))^{-1}$$
satisfies the requirements.
\end{proof}

\subsection{Markov measures and the Radon-Nikodym problem.}

We assume that the space $X$ is compact. The above theorem reduces the Radon-Nikodym problem
for a continuous cocycle
$c:R\rightarrow{\bf R}$, where $R$ is an AF equivalence relation, to the case when
$c$ is a quasi-product cocycle. Thus we assume that $R$ is the tail equivalence relation
on the path space of a Bratteli diagram $(V,E)$ and that $D\in Z^1(R,{\bf R}_+^*)$
is the quasi-product cocycle defined by the labelling $\Phi: E\rightarrow {\bf R}_+^*$.
Then we have the following elementary result:

\begin{prop} The solutions of the (RN) equation $D_\mu=D$ are the Markov measures
$$\mu(Z(x_1\ldots x_n))=\rho_0(s(x_1))p_1(x_1)\ldots p_n(x_n)$$
where $p_n(e)=(\rho_{n-1}\circ s(e))^{-1}\Phi_n(e)\rho_n\circ r(e)$ for
$n=1,2,\ldots, e\in E(n)$ and $\rho:V\rightarrow {\bf R}_+^*$ is a normalized
$\Phi$-harmonic function, i.e. a solution of
\begin{eqnarray}
\rho_{n-1}(v)&=&\sum_{s(e)=v}\Phi(e)\rho_n\circ r(e) \quad (n=1,2,\ldots\ v\in
V(n-1))\\ 
1&=&\sum_{V(0)}\rho_0(v)
\end{eqnarray}
\end{prop}
\begin{proof}
We first check that the Markov measure $\mu$ defined by the transition
probability $p$ and the initial measure $\mu_0$ is quasi-invariant
with Radon-Nikodym
derivative
$D$ given by
$$D(az,bz)={\mu_0(s(a))p(a)\over\mu_0(s(b)p(b)},$$
where $a,b$ are finite paths and $p(a_1\ldots a_n)= p_1(a_1)\ldots p_n(a_n)$.
This means that the equality
$$\int_{r(S)} f(x) d\mu(x)=\int_{s(S)} f\circ r(Sy) D(Sy) d\mu(y)$$
holds for every open bissection $S$ and every
$f\in C_c(X)$ with support contained in $r(S)$. It suffices to consider
finite changes of coordinates
$$S=S(a,b):bx\in Z(b)\mapsto ax\in Z(a)$$
where $a,b$ are finite paths in the Bratteli diagram starting from the level $0$
and ending at the same vertex $v$ and $f=\chi_{Z(ac)}$, where $c$ is a finite path
starting at $v$. Then the left handside is $\mu(Z(ac))$ while the right handside is
$${\mu_0(s(a))p(a)\over\mu_0(s(b)p(b)}\mu(Z(bc))=\mu(Z(ac)).$$
Conversely, suppose that $\mu$ is a quasi-invariant probability measure having $D$
as Radon-Nikodym derivative. The above argument gives the equality
$${\mu(Z(ac))\over \Phi(a)}={\mu(Z(bc))\over \Phi(b)}.$$
In particular $\mu(Z(a))/ \Phi(a)$ depends only on its end vertex $v=r(a)$.
Therefore there exists $\rho: V\rightarrow{\bf R}_+^*$ such that
$\mu(Z(a))=\Phi(a)\rho(r(a))$. The condition $(1)$ results from the additivity of
$\mu$:
$$\mu(Z(a))=\sum_{s(e)=r(a)}\mu(Z(ae))$$
while the condition $(2)$ holds because $\mu$ is a probability measure.
\end{proof}

For example, let us consider the stationary Bratteli diagram associated with the
full one-sided shift: for $n\in{\bf N}$, $V(n)=\{1,\ldots,d\}$ and for $n\ge 1$,
$E(n)=\{(v,w): v,w\in \{1,\ldots,d\}$. A matrix $A\in M_d({\bf R}_+^*)$ defines a
stationary $\Phi: E\rightarrow {\bf R}_+^*$ such that $\Phi_n(v,w)=A_{v,w}$. The
above equation becomes:
\begin{eqnarray}
\rho_{n-1}&=&A\rho_n\quad (n=1,2,\ldots)\nonumber \\ 
1&=&\sum_1^d\rho_0(v)\nonumber
\end{eqnarray}
It admits a unique solution $\rho_n=\lambda^{-n}\rho_0$, where $\lambda$ is the
Perron-Frobenius eigenvalue of $A$ and $\rho_0$ is the normalized Perron-Frobenius
eigenvector.

\section{Stationary cocycles.}

We return now to the situation described in the introduction. We may consider
the case of a local homeomorphism $T$ of a compact metric space $X$ onto
itself rather than the special case of the one-sided full shift.

\subsection{Singly generated dynamical systems}
Given $(X,T)$ as above, we define as in the introduction:
\begin{eqnarray}
G(X,T)&=&\{(x,m-n,y)\in X\times{\bf Z}\times X:m,n\in{\bf N},
T^mx=T^ny\}\nonumber\\ 
R(X,T)&=&\{(x,y)\in X\times X:\exists n\in{\bf N},
T^nx=T^ny\}\nonumber.
\end{eqnarray}
They have a natural structure of \'etale locally compact groupoid (see for example
\cite{dea:groupoids} or \cite{ren:cuntzlike}). In fact $R(X,T)=\cup R_n$, where
$R_n$ is obtained by fixing $n$, is an AP equivalence relation. Moreover, $G(X,T)$
is the semi-direct product of $R(X,T)$ by the endomorphism induced by $T$. Given a
group
$A$, there is a one-to-one correspondence between functions $\varphi:X\rightarrow
A$ and one-cocycles $c:G(X,T)\rightarrow A$. It is given by
$$c_\varphi(x,m-n,y)=\varphi(x)\varphi(Tx)\ldots\varphi(T^{m-1}x)
\varphi(T^{n-1}y)^{-1}\ldots\varphi(Ty)^{-1}\varphi(y)^{-1}.$$
Moreover, if $A$ is a locally compact group, $c_\varphi$ is continuous iff
$\varphi$ is continuous. We shall denote by $Z^1_{cont}(G,A)$ the space of
continuous cocycles from $G$ to $A$. Here are some elementary observations (for
brevity, we write $G=G(X,T), R=R(X,T)$):
\begin{enumerate}
\item $c\in Z^1_{cont}(R,{\bf R})$ extends to $G$ iff there exists $\varphi\in
C(X)$ such that
$$c(x,y)-c(Tx,Ty)=\varphi(x)-\varphi(y).$$
Then $c={c_\varphi}_{|R}$.
\item ${c_\varphi}_{|R}={c_\phi}_{|R}$ iff $\varphi-\phi$ is invariant under $R$.
(If $R$ has a dense orbit, $\varphi-\phi$ must then be constant).
\item If $\mu$ is a measure on $X$ quasi-invariant under $G$, it must be
quasi-invariant under $R$. If $D_\mu$ is its Radon-Nikodym-derivative with respect to $G$,
then its Radon-Nikodym-derivative with respect to $R$ is the restriction
${D_\mu}_{|R}$.
\item If the Radon-Nikodym problem for $D\in Z^1_{cont}(R,{\bf R}_+^*)$ has a unique solution
and $D$ extends to $G$, then $\mu$ is quasi-invariant under $G$. 
\end{enumerate}
A cocycle $c\in Z^1_{cont}(R,A)$ which extends to $G(X,T)$ will be called {\it
stationary}.
\vskip5mm

These facts suggest the following scheme to solve the Radon-Nikodym problem for $G(X,T)$ and
the cocycle $e^{c_\varphi}\in Z^1_{cont}(G(X,T),{\bf R}_+^*)$, where $\varphi\in
C(X,{\bf R})$.

$(a)$ Solve the Radon-Nikodym problem for $R(X,T)$ and
the stationary cocycle  $e^{{c_\varphi}_{|R}}\in Z^1_{cont}(R(X,T),{\bf R}_+^*)$.

$(b)$ Check the quasi-invariance under $G(X,T)$ of the solutions $\mu$
found in the first step (this is automatic when $a$ has a unique solution). If this
is the case, compare
$D_\mu$ and
$e^{c_\varphi}$.

The advantage of this scheme is to reduce the problem to an AP (or AF)
equivalence relation, for which approximation techniques matingale convergence
theorems are available. If the problem
$a$ admits a unique solution, then there is a unique constant $\beta\in{\bf R}$,
depending on $\varphi$ such that $\varphi-\beta$ defines a Radon-Nikodym cocycle for $G(X,T)$.
This constant $\beta$ is the topological pressure $p(\varphi)$ of $\varphi\in
C(X)$. Instead of developing here this approach, we will apply directly the
Perron-Frobenius-Ruelle theorem which solves the Radon-Nikodym problem under suitable
assumptions. 

\subsection{The Perron-Frobenius-Ruelle theorem.}

The Perron-Frobenius-Ruelle theorem gives directly a solution of the Radon-Nikodym
problem for $G(X,T)$ rather than for $R(X,T)$. It requires suitable assumptions
both on the dynamical system
$(X,T)$ and the function
$\varphi\in C(X,{\bf R})$.

In the literature, various assumptions (or may be different expressions of
these assumptions) are made. The version of the Perron-Frobenius-Ruelle theorem
given by P. Walters in \cite{wal:equilibrium} is particularly well-adapted to
our setting . His assumptions on the dynamical system
$(X,T)$ are in fact weaker than what we need here. We assume that:
\begin{enumerate}
\item $T$ is expansive: there exists $\epsilon>0$ such that, if $x\not=y$, there
exists $n\in{\bf N}$ with $d(T^nx,T^ny)>\epsilon$;
\item $R(X,T)$ is minimal: all the orbits with respect to $R(X,T)$ ($x\sim y$ iff
there exist $n\in{\bf N}$ such that $T^nx=T^ny$) are dense.
\end{enumerate}

The assumption on the cocycle is that $\sum_1^\infty var_n(\varphi)<\infty$,
where
$$var_n(\varphi)=\sup\{|\varphi(x)-\varphi(y)|: d_n(x,y)<1\}$$
and
$$d_n(x,y)=d(x,y)+d(Tx,Ty)+\ldots +d(T^{n-1}x,T^{n-1}y).$$

The basic tool in this theory is the Ruelle operator ${\mathcal
L}_\varphi:C(X)\rightarrow C(X)$ defined by
$${\mathcal L}_\varphi f(x)=\sum_{Tx=y}e^\varphi(y)f(y) .$$

\begin{thm}(Perron-Frobenius-Ruelle, \cite{wal:equilibrium}). Under
the above assumptions, the eigenvalue problem 
$$ ^t{\mathcal L}_\varphi\mu=\lambda\mu , $$
where $\lambda>0$ and $\mu$ is a probability measure, admits one and only one
solution $(\lambda,\mu)$. Moreover,  $\log\lambda$ is the
topological pressure $p(\varphi)$ of $\varphi$.
\end{thm}

The following result expresses
the Radon-Nikodym problem in a form similar to that given in Proposition 3.2. for
quasi-product cocycles. Again, it is an elementary result, which is well known  in
the Perron-Frobenius-Ruelle theory.

\begin{prop} The solutions of the ($RN$) equation $D_\mu=e^{c_\varphi}$ are the
the solutions of
$$ ^t{\mathcal L}_\varphi\mu=\mu.$$
\end{prop}
\begin{proof}
Suppose that $\mu$ is a quasi-invariant probability measure admitting
$D=e^\varphi$ as its Radon-Nikodym derivative. Let $J(x)=e^{\varphi(x)}=D(x,1,Tx)$ be its
Jacobian. We have for all $f\in C(X)$:
\begin{eqnarray}
<f,^t{\cal L}_\varphi \mu>&=&<{\cal L}_\varphi f,\mu>\nonumber\\
&=& \int_X\sum_{Ty=x}J(y)f(y) d\mu(x)\nonumber\\
&=&\int g d s^*\mu\nonumber
\end{eqnarray}
where $g=((J f)\circ r)\chi_S$ and $S=\{(x,1,Tx): x\in X\}$. By
quasi-invariance of $\mu$, this is
\begin{eqnarray}
\int g D^{-1} dr^*\mu &=& \int_X J(x)f(x)J(x)^{-1} d\mu(x)\nonumber\\
&=& <f,\mu>\nonumber
\end{eqnarray}
Suppose now that the probability measure $\mu$ satisfies $ ^t{\mathcal
L}_\varphi\mu=\mu $. Because of the equality
$${\cal L}_\varphi (f\circ T)=({\cal L}_\varphi 1)f,$$
we have $T_*\mu=({\cal L}_\varphi 1)\mu$, where $T_*\mu(A)=\mu(T^{-1}A)$. In
particular $\mu(A)=0\Rightarrow \mu(T^{-1}A)=0$. If $\mu(A)=0$, we also have
$$\int_X\sum_{Ty=x}e^{\varphi(y)}\chi_A(y) d\mu(x)=0$$
and this implies $\mu(TA)=0$. These two properties imply the quasi-invariance of
$\mu$ under $G(X,T)$. Let $D$ be its Radon-Nikodym derivative and $J(x)=D(x,1,Tx)$ be its
Jacobian. As above, we have for all $f\in C(X)$,
\begin{eqnarray}
<f,\mu>&=& <f,^t{\cal L}_\varphi \mu>\nonumber\\
&=& \int_X e^{\varphi(x)}f(x)J(x)^{-1} d\mu(x)\nonumber
\end{eqnarray}
which implies that $J(x)=e^{\varphi(x)}$ for $\mu$-a.e. $x$ and that
$D=e^{c_\varphi}$.
\end{proof}

\begin{cor} Under the assumptions of the PFR theorem, the ($RN$) equation
$D_\mu=e^{c_\varphi}$ admits a solution iff the topological pressure $p(\varphi)$
of
$\varphi$ is zero. In that case, the solution is unique.
\end{cor}

Note that, under the above assumptions on $(X,T)$ and $\varphi$, the set 
$\{D_{\varphi+\beta} :\beta\in{\bf R}\}$ contains one
and only one Radon-Nikodym cocycle. It corresponds to the value $\beta=-p(\varphi)$. However,
this does not quite solve the usual KMS problem, which concerns the set
$\{D_{-\beta\varphi} :\beta\in{\bf R}\}$. It remains to see if the equation
$p(-\beta\varphi)=0$ has a unique solution in $\beta$. This is done in particular
cases in \cite{eva:on} and \cite{exe:RPF} under the assumption that $\varphi$ is
strictly positive. This can be done in a more general framework by using the
convexity properties of the pressure function. Note that in the case of the usual
gauge action,
$\varphi\equiv 1$ and the equation $p(-\beta\varphi)=0$ has the unique solution
$\beta=p(0)$.


\begin{thebibliography}{10}

\bibitem{aeg:amenable} S.~Adams, G.A.~Elliott and T.~Giordano : {\it Amenable
actions of  groups}, Trans. Amer. Math. Soc., {\bf 344} (1994), 803-822.

\bibitem{av:endo} V.~Arzumanian and A.~Vershik: {\it Star-algebras
associated with endomorphisms} in {\it Operator algebras and group
representations}, Proc.~Int.~Conf., Vol.~{\bf I}, Pitman Boston, 1984, 17-27.

\bibitem{ar:examples} V.~ Arzumanian and J.~Renault: {\it Examples of pseudogroups
and their $C^*$-algebras}, in Operator Algebras and Quantum Field Theory, 
S.~Doplicher, R.~Longo, J.~E.~Roberts and L.~Zsido, editors, International Press
1997, 93-104.

\bibitem{cfw:amenable}  A.~Connes, J.~Feldman and B.~Weiss: {\it An amenable
equivalence relation is generated by a single transformation}, J.~Ergodic Theory and
Dynamical Systems, {\bf 1} (1981), 431-450.

\bibitem{cun:simple}  J.~Cuntz: {\it Simple $C^*$-algebras generated by
 isometries}, Comm. Math. Phys., {\bf  57} (1977), 173-185.

\bibitem{ckr:markovI}  J.~Cuntz and W.~Krieger: {\it A class of $C^*$-algebras and
topological Markov chains},
Inventiones Math., {\bf  56} (1980), 251-268.

\bibitem{dea:groupoids} V.~Deaconu: {\it Groupoids associated with
endomorphisms}, Trans. Amer. Math. Soc., {\bf 347} (1995), 1779-1786.

\bibitem{dye:I} H.~Dye: {\it On groups of measure preserving transformations I},
 Amer. J. Math., {\bf 81} (1959), 119-159.

\bibitem{dye:II} H.~Dye: {\it On groups of measure preserving transformations II},
 Amer. J. Math., {\bf 85} (1963), 551-576.

\bibitem{eg:amenable} G.A.~Elliott and T.~Giordano: {\it Amenable
actions of discrete groups}, Ergodic Theory Dynam. Systems, {\bf 13} (1993),
no 2, 803-822.

\bibitem{efw:kms} M.~Enomoto, M.~Fujii and Y.~Watatani: {\it KMS states for
gauge action on $O_A$}, Math. Japon., {\bf 29} (1984), no 4,
607-619.

\bibitem{eva:on}  D.~Evans: {\it On $O_n$}, Publ. RIMS Kyoto Univ. {\bf 16} 
(1980), 915-927.

\bibitem{eva:on+}  D.~Evans: {\it On $O_{n+1}$}, Invent.~math. {\bf 93}  (1989),
1-13.

\bibitem{eva:markov}  D.~Evans: {\it The $C^*$-algebras of topological Markov
chains}, Invent.~math. {\bf 93}  (1989), 1-13.

\bibitem{eva:quasiproduct}  D.~Evans: {\it Quasi-product states on C$^*$-algebras},
Invent.~math. {\bf 93}  (1989), 1-13.

\bibitem{exe:RPF} R.~Exel:{\it KMS states for generalized gauge actions on
Cuntz-Krieger algebras}, Preprint, (October 2001).

\bibitem{el:infinitematrices}  R.~Exel and M.~Laca: {\it Cuntz-Krieger algebras
for infinite matrices}, J.reine angew. Math.
{\bf 512} (1999),119-172.

\bibitem{el:KMS} R.~Exel and M.~Laca: {\it Partial dynamical systems and the KMS
condition}, preprint 2000.

\bibitem{fm:relations}  J.~Feldman and C.~Moore: {\it Ergodic equivalence relations,
cohomologies, von Neumann algebras, I and II},
Trans. Amer. Math. Soc., {\bf 234}, no. 2 (1977), 289-359.

\bibitem{gps:orbitequivalence}  T.~Giordano, I.~Putnam and C.~Skau: {\it
Topological orbit equivalence and $C^*$-crossed products}, J.reine angew. Math.
{\bf 469} (1995),51-111.

\bibitem{kv:walk}  V.~Kaimanovich and A.~Vershik: {\it Random walks on discrete
groups:boundary and entropy}, Ann. Probab., {\bf 11} (1983), no 3, 457-490.

\bibitem{kri:flows}  W.~Krieger: {\it On ergodic flows and the isomorphism of
factors}, Math. Ann., {\bf  223} (1976), 19-70.

\bibitem{op:KMS}  D.~Olesen and G.~Pedersen: {\it  Some $C^*$-dynamical systems
with a single KMS state}, Math. Scand. {\bf 42}
(1978), 179-197.

\bibitem{ren:approach}  J.~Renault: {\it A groupoid approach to
$C^*$-algebras}, Lecture Notes in Mathematics, Vol.~{\bf 793}
Springer-Verlag Berlin, Heidelberg, New York (1980).

\bibitem{ren:cuntzlike}  J.~Renault: {\it  Cuntzlike-algebras}, J.~Operator Theory
{\bf 25} (1991), 3-36.

\bibitem{rue:hyperbolic}  D.~Ruelle: {\it Noncommutative algebras for hyperbolic
diffeomorphisms}, Inventiones Math. {\bf 93}  (1989), 1-13.

\bibitem{sch:cocycles} K.~Schmidt: {\it Cocycles on ergodic transformation
groups}, MacMillan (Company of India, Ltd) Delhi, 1977.

\bibitem{wal:equilibrium} P.~Walters: {\it Invariant measures and equilibrium
states for some mappings which expand distances} Trans. Amer. Math. Soc., {\bf
236}, 121-153.

\bibitem{zim:poisson} R.~Zimmer : {\it Amenable ergodic group actions and an
application to Poisson boundaries of random walks}, J. Funct. Anal., 
{\bf 27} (1978), 350-372.





\end{thebibliography}
\end{document}